\documentclass[12pt]{amsart}

\title{Functions with image in a strip}
\author[T.  Bhattacharyya, A. G. O'Farrell, \\ S. Rastogi and Vijaya Kumar U.]{Tirthankar Bhattacharyya, Anthony G. O'Farrell, \\ Shubham Rastogi and Vijaya Kumar U.}

\newcommand{\Addresses}{ \bigskip
  \footnotesize
  
  T. ~Bhattacharyya and Vijaya Kumar U., \textsc{Department of Mathematics, Indian Institute of Science, Bangalore 560012.}  \par\nopagebreak \textit{E-mail addresses}:  \texttt{tirtha@iisc.ac.in} and \texttt{vijayak@iisc.ac.in}

\medskip

A. G. ~O'Farrell, \textsc{Department of Mathematics and Statistics, Maynooth University, Naas, Co. Kildare, W23HW31 Ireland.} \par\nopagebreak \textit{E-mail address}: \texttt{anthony.ofarrell@mu.ie}

\medskip

S. ~Rastogi, \textsc{Department of Mathematics, Indian Institute of Technology Bombay, Mumbai 400076, India.}  \par\nopagebreak \textit{E-mail address}: 
\texttt{shubhamr@math.iitb.ac.in}}

\RequirePackage{datetime2,graphicx}
\usepackage[mathscr]{eucal}
\date{\DTMnow}
\newcommand\C{\mathbb{C}}
\newcommand\D{\mathbb{D}}
\newcommand\R{\mathbb{R}}

\newcommand{\z}{\mathbf{z}}

\newcommand\E{\mathcal{E}}

\newcommand\Po{\mathcal{P}}

\newcommand{\V}{\mathscr V}
\renewcommand{\S}{\mathscr S}

\newcommand\ignore[1]{}
\newtheorem{theorem}{Theorem}
\newtheorem{lemma}{Lemma}

\newtheorem*{conjecture*}{Conjecture}
\newcommand\be{\begin{equation}}
\newcommand\ee{\end{equation}}
\begin{document}

\begin{abstract}
	We consider special functions.
\end{abstract}

\maketitle

\vspace*{5mm}

\begin{centerline}
{Dedicated to the memory of K. R. Parthasarathy}
\end{centerline}
\begin{centerline}
{whose teaching has had a great influence on the first named author.}
\end{centerline}

\section{Introduction}\label{S:1}
\subsection{The question}
Let $\Re$ and $\Im$ denote respectively the real part and the imaginary part of a complex number. We begin with the following question. Do the conditions $0\le \Re F(z)\le 1$ and $F(z)+z\overline{F(z)}$ holomorphic oblige a function $F:\D\to\C$ to be constant?

For $F:\D\to\C$, we define 
$$ (LF)(z):= F(z)+z\overline{F(z)}, \ \forall z\in\D.$$
This is an $\R$-linear function on $\C^{\D}$, which preserves continuity, differentiability of any order, and real analyticity.

Consider the strip
$$ \text{ Strip } := \{z\in\C: 0\le\Re z\le 1\}.$$

\begin{theorem}\label{T:1}
If $F:\D\to\text{ Strip } $ and $LF$ is holomorphic
then $F$ is constant.
\end{theorem}

In Section \ref{S1:convexity} we prove this by appealing to 
known facts about the extreme points of a certain convex set of functions.

When one considers holomorphic functions with positive real part, it is not surprising that the representation by Herglotz plays a role. Using that we prove the following theorem in Section \ref{S2:Herglotz}.

\begin{theorem}\label{T:2}
Let $\mathcal H$ be a complex, separable Hilbert space and let $F:\D\to\mathcal B (\mathcal H)$ be a function satisfying
\begin{enumerate}
\item the real part of $F$ is a positive contraction,
\item the function $F(z) + z F(z)^*$ is holomorphic.
\end{enumerate} 
Then $F$ is constant.
\end{theorem}

We then apply Theorem \ref{T:2} to find all factorizations of the shift semigroup.

\section{The convexity argument}\label{S1:convexity}
Let $F=b+ia$, where $b:\D\to [0,1]$ and $a:\D\to\R$, and
assume $g=LF$ is holomorphic on $\D$.

Then 
$$ g(z) = (1+z)b(z) + i(1-z) a(z). $$
Let
$$ \phi(z):= \frac{1+z}{1-z}. $$
Define 
$$h_1(z):= b(z)\phi(z) + ia(z),\ h_2(z):=\phi(z)-h_1(z). $$
Then
$$ h_1(z) = \frac{g(z)}{1-z} $$
is holomorphic in $\D$ and
$$ \Re h_1(z) = b(z) \Re\phi(z) = b(z)\frac{1-|z|^2}{|1-z|^2} \ge0 $$
on $\D$.

If $h_1$ is nonconstant, then $\Re h_1(z)>0$ on $\D$.
Since $z \mapsto \frac{1+z}{1-z}$ maps 
the unit disc conformally onto 
the right half-plane, 
there exists a 
holomorphic function
$\psi_1:\D\to\D$ such that 
$$ h_1(z) = \frac{1+\psi_1(z)}{1-\psi_1(z)}. $$
If, on the other hand, $h_1$ is constant on $\D$,
then there is a constant $\psi_1$ in the closed unit disk
such that 
$$ h_1(z) = \frac{1+\psi_1}{1-\psi_1}. $$
So in all cases we can say that 
there exists a 
holomorphic function
$\psi_1$ on $\mathbb D$ such that $|\psi_1|\le 1$ and 
$$ h_1(z) = \frac{1+\psi_1(z)}{1-\psi_1(z)}. $$

We also have  $\Re h_2\ge 0$ on $\D$, so in the same way
there exists a 
holomorphic function
$\psi_2$ on $\mathbb D$ such that $|\psi_2|\le 1$ and 
$$ h_2(z) = \frac{1+\psi_2(z)}{1-\psi_2(z)}. $$

Thus
\begin{equation}\label{E:cc}
\frac{1+\psi_1(z)}{1-\psi_1(z)}+
\frac{1+\psi_2(z)}{1-\psi_2(z)}
 = \frac{1+z}{1-z}.
\end{equation} 

Let
$\Po$ denote the set of holomorphic functions $f$ on $\D$
that have $f(0)=1$ and $\Re f>0$ on $\D$. 
Let
$$ f_j(z) = \frac1{\Re h_j(0)}
\left[ h_j(z)-i\Im h_j(0) \right].
$$
Then $f_1,f_2\in\Po$ and 
$$ (\Re h_1(0)) \cdot f_1(z) 
+ (\Re h_2(0)) \cdot f_2(z) 
= \phi(z). $$
This expresses $\phi(z)$ as a convex combination of
$f_1$ and $f_2$. But by \cite[Theorem 1.5]{Schober},
$\phi$ is an extreme point of the convex set $\Po$,
so we conclude that $f_1=f_2=\phi$.
So
$$\frac{1+\psi_j(z)}{1-\psi_j(z)}
= 
i \Im \frac{1+\psi_j(0)}{1-\psi_j(0)}
+
\Re\left(\frac{1+\psi_j(0)}{1-\psi_j(0)}\right)
\frac{1+z}{1-z}, $$
So, with $h:=h_1$,
\be\label{E:7}
h(z)=\Re h(0) \phi(z) + i\Im h(0), $$
$$ \Re h(z)=\Re h(0)\Re\phi(z) 
=\Re h(0) \frac{1-|z|^2}{|1-z|^2}.
\ee
Since $h(z)=b(z)\phi(z)+ia(z)$, we have
\be\label{E:8}
\Re h(z) = b(z)\frac{1-|z|^2}{|1-z|^2}.
\ee
By \eqref{E:7} and \eqref{E:8},
$$ b(z)\frac{1-|z|^2}{|1-z|^2}= \Re h(0)\frac{1-|z|^2}{|1-z|^2},\ \forall z\in\D.
$$
So $b(z)$ is the constant 
$$ \Re h(0) = 
\Re \frac{1+\psi_1(0)}{1-\psi_1(0)}.
$$
Now $h(z)$ is analytic and $b\phi(z)$ is analytic, so
$h(z)-b\phi(z)$ is analytic.  But this means that $ia(z)$ is
analytic and imaginary-valued.  So $a(z)$ is constant.

This proves Theorem \ref{T:1}.

\section{The operator valued case - the Herglotz argument} \label{S2:Herglotz}

The proof of Theorem \ref{T:2} requires the Herglotz representation of a holomorphic function on the disc with positive real part. We quote from \cite{KBS}.

\begin{theorem}[Theorem 2.1 in \cite{KBS}]
A function $A:\mathbb{D} \longrightarrow \mathcal{B}(\mathcal{H})$ is analytic with non-positive real part for all $z \in \mathbb{D}$ if and only if it admits a representation:
\begin{equation}
A(z)=i\textup{Im}A(0)+ \int_{\mathbb T}\frac{z+\exp(it)}{z-\exp(it)}d\Sigma(\exp(it))\label{Az}\\
\end{equation}
for $z \in \mathbb{D}$, with a positive operator valued measure $\Sigma$. Furthermore, such a $\Sigma$ is unique.
\end{theorem}

As before, assume that $F(z) = B(z) + iA(z)$ where the functions $B$ and $A$ defined on $\mathbb D$ take values in the set of self-adjoint elements of $\mathcal B (\mathcal H)$. Then, $F(z) + z F(z)^* = (1 + z)B(z) + i(1 - z)A(z)$. Call this $g(z)$. Since $g(z)$ is holomorphic on $\mathbb D$, so is $h_1(z) = \phi(z) B(z) + iA(z)$ whose real part is $B(z) \frac{1 - |z|^2}{|1-z|^2}$.  Let $h_2(z) = \phi(z) I - h_1(z)$ whose real part is $(I - B(z)) \frac{1 - |z|^2}{|1-z|^2}$. Since $B(z)$ is a positive contraction, both these real parts are positive operators for each $z \in \mathbb D$. 

Hence, we have measures $\Sigma_1$ and $\Sigma_2$ taking values in the cone of positive elements of $\mathcal B(\mathcal H)$ such that 
$$ h_j(z) = i\textup{Im}h_j(0)+ \int_{\mathbb T}\frac{z+\exp(it)}{z-\exp(it)}d\Sigma_j(\exp(it)), j=1,2.$$
Adding them, we get 
$$\varphi(z) =  \int_{\mathbb T}\frac{z+\exp(it)}{z-\exp(it)}d\Sigma(\exp(it))$$
where $\Sigma = \Sigma_1 + \Sigma_2$. Now, we invoke the uniqueness of the measure in Herglotz representation. That implies that $\Sigma$ must be the Dirac measure with mass at $1$. Since $\Sigma_1$ and $\Sigma_2$ are non-negative operator valued measures, they are concentrated at $\{1\}$. Let $B_j = \Sigma_j(\{1\})$ and $A_j = \textup{Im}h_j(0)$. Then,
$$ B_1 + B_2 = I \text{ and } h_j(z) = i A_j + \frac{1+z}{1-z}B_j$$
for $j=1,2$. 

Now we prove a small lemma which will be required to complete the proof. 

\begin{lemma}
If $F,g:\D\to\mathcal B (\mathcal H)$, then
$$ F(z) + zF(z)^*=g(z) \text{ if and only if } (1-|z|^2)F(z)=g(z)-zg(z)^*. $$
\end{lemma}
\begin{proof} If $F(z) + zF(z)^*=g(z)$, then
\begin{align*}
F(z)^* + \bar z F(z) &= g(z)^*,\\
zF(z)^* + |z|^2 F(z) &= z g(z)^*
\end{align*}
and subtracting, we get $(1-|z|^2) F(z) = g(z)-z g(z)^*$.
The converse is routine. 
\end{proof}

Resuming the proof of Theorem \ref{T:2}, we note that 
$$h_1(z) = i A_1 + \frac{1+z}{1-z}B_1$$ 
and hence 
$$g(z) = (1+z) B_1 +i(1-z)A_1 = B_1 + iA_1 + z(B_1 - iA_1) = C + zC^*$$
where $C = B_1 + iA_1$. Consequently,
$$F(z)=\frac{g(z)-zg(z)^*}{1-|z|^2} = C.$$

\section{Application}

As an application, we deduce a known result about factorizations.  Let $\mathcal E$ be a finite-dimensional complex separable Hilbert space. The shift semigroup $\S^\E = (S_t^\E)_{t\ge 0}$ on $L^2(\R_+,\E)$ is defined by
	\[(S_t^\E f)x=\begin{cases}
		f(x-t) &\text{if } x\ge t,\\
		0 & \text{otherwise}
	\end{cases}\]
	for $t\ge 0.$ When $\E=\C$ we just write $\S= (S_t)_{t\ge 0}$ instead of $\S^\E.$ We shall change the underlying space by a unitary conjugation so that we can use the theory of holomorphic functions.  If $S^\E$ is the co-generator, i.e., the Cayley transform of the generator of the semigroup  $\S^\E$, then we define a unitary operator    $W_\E:L^2(\R_+,\E)\to H^2_\D(\E)$    by
	\begin{equation}\label{shift}
		W_\E({(S^\E)}^n(\sqrt{2}e^{-x}f))=z^nf\text{ for } f\in \E, n\ge 0.
	\end{equation}
Let $\varphi_t(z) = \exp (-t\phi(z))$. Then it is known that	
\begin{equation}\label{conj}
		W_\E S_t^\E W_\E^*=M_{\varphi_t\circ \z^\E} \text{ for all } t\ge 0,
	\end{equation}
where $\z^\E = z I_{\E}$. 

A pair of isometric semigroups $(\V_1, \V_2)$ is said to be a {\em factorization} of an isometric semigroup $\V$ if $\V_1$ and $\V_2$ are commuting and
	\begin{equation}\label{defn:fact}
		V_t = V_{1,t}V_{2,t}
		\text{ for all }t\ge 0.
	\end{equation} 
By using the unitary operator $W_\E$, a factorization of the shift semigroup translates to finding a pair of $\mathcal B (\E)$-valued bounded holomorphic functions $\psi_1$ and $\psi_2$ on the unit disc $\D$ with  $$\|\psi_j\|_\infty=\sup \{\|\psi_j(z)\|:z\in \D\}\le 1$$ for $j=1,2$ and  satisfying any one of the two equivalent statements:

(a)  the semigroups $(M_{\varphi_t\circ \psi_1})_{t\ge 0}$ and $(M_{\varphi_t\circ \psi_2})_{t\ge 0}$ factorize $(M_{\varphi_t\circ \z^\E})_{t\ge 0}$,

(b) $1\notin \sigma_p(\psi_j(z))$ for any $z\in \D$ and  $\psi_1$ and $\psi_2$  satisfy
		\begin{equation}\label{eq:master}
			(I+\psi_1(z))(I-\psi_1(z))^{-1}+(I+\psi_2(z))(I - \psi_2(z))^{-1}=\phi(z)I
		\end{equation}

		for all $z\in \D.$ Compare this with \eqref{E:cc}.

We call such a pair $(\psi_1, \psi_2)$ of functions a factorizing pair. Proving the above needs some effort and it will be a repetition if we went into the details of the proofs here. The details are available in \cite{BRSU}. 

As a consequence of \eqref{eq:master}, we have
$$ 0 \le \Re (I+\psi_j(z))(I-\psi_j(z))^{-1} \le \Re \phi(z)I_\E$$
for any factorizing pair $(\psi_1, \psi_2)$. Let
$$ B(z) = \Re (I+\psi_1(z))(I-\psi_1(z))^{-1} / \Re \phi(z). $$
Clearly, $B(z)$ is a positive contraction for each $z$ in $\D$. Note that
$$ (I+\psi_1(z))(I-\psi_1(z))^{-1} =  \phi(z)B(z) + i A(z)$$
where $A(z) = \Im (I+\psi_1(z))(I-\psi_1(z))^{-1} - \Im (\phi(z))B(z)$ is a function on $\D$ taking values in self-adjoint operators on $\mathcal E$. Since 
$$(I+\psi_1(z))(I-\psi_1(z))^{-1}$$ 
is holomorphic, the stage is set for application of Theorem  \ref{T:2}. Thus, we get that $A(z) =- A$ for all $z$ in $\D$ and $B(z) = B$ for all $z$ in $\D$ where $A$ is a self-adjoint operator and $B$ is a positive contraction. In particular, we have 
 \begin{align*}
 	\varphi_t\circ\psi_1(z)&=e^{t[iA+(z+1)(z-1)^{-1}B]} \text{ and}\\
 	\varphi_t\circ\psi_2(z)&=e^{t[-iA+(z+1)(z-1)^{-1}(I-B)]}
 \end{align*} 
 	for all $z\in \D.$ 
 We have proved the following theorem which is a special case of the main result of \cite{BRSU}.

\begin{theorem}Let  $(\V_1,\V_2)$ be a factorization of  $(M_{\varphi_t\circ\z^\E})_{t\ge 0}.$ Then there exist  self-adjoint matrices $A$ and $B$ in $\mathscr B(\E)$ with $0\le B\le I$ such that  
\begin{equation}\label{eq:Vjt}
V_{j,t}=M_{\varphi_{j,t}} \text{ for all }t\ge 0, j=1,2,
\end{equation} 
where $\varphi_{j,t}\in H^\infty_\D(\mathscr B(\E))$ is defined by 
\begin{equation}\label{eq:phi-j}\varphi_{1,t}(z)=e^{t[iA+(z+1)(z-1)^{-1}B]}\text{ and }\varphi_{2,t}(z)=e^{t[-iA+(z+1)(z-1)^{-1}(I-B)]}.
\end{equation}

Conversely, suppose $A,B\in \mathscr B(\mathcal H)$ are self-adjoint and $0\le B\le I,$ if we define $\varphi_{j,t}$ as in \eqref{eq:phi-j} then $\varphi_{j,t}\in H^\infty_\D(\mathscr B(\E))$ for $j=1,2,t\ge 0$ and  $(\V_1,\V_2)$ as defined in \eqref{eq:Vjt} is a factorization of  $(M_{\varphi_t\circ\z^\E})_{t\ge 0}.$
\end{theorem}

	\noindent{\bf Acknowledgements}
	The first named author gratefully acknowledges the hospitality of the Maynooth University and the J C Bose fellowship JCB/2021/000041 of SERB. The fourth named author is supported by the  D S Kothari postdoctoral fellowship MA/20-21/0047. Research is supported by the DST FIST program - 2021 [TPN - 700661]. All authors thank the referee for very useful suggestions for future work.

\Addresses

\end{document}